\documentclass[ams,12pt,twoside]{article}

\usepackage[tbtags]{amsmath}
\usepackage{amssymb}
\usepackage{amsfonts}
\usepackage{euscript}
\usepackage{longtable}
\usepackage[dvips]{graphicx}
\usepackage{multicol}
\usepackage{color}
\usepackage{rotating}
\usepackage{hhline}
\usepackage[matrix,arrow,curve]{xy}

\makeglossary
\makeindex

\newcommand{\bs}{\boldsymbol}
\newcommand{\mf}{\mathfrak}
\newcommand{\mr}{\mathrm }
\newcommand{\mt}{\mathrm{Tr} }

\newcommand{\ld}{\ldots }

\newcommand{\mre}{\mathrm{e} }

\begin{document}

\title{\bf A method for construction  \\ of Lie group invariants}

\author{Yu.~Palii \\
\small{\it Laboratory of Information Technologies,} \\
\small{\it Joint Institute for Nuclear Research, Dubna, Russia} \\
\small{\it and} \\
\small{\it {Institute of Applied Physics, Chisinau, Moldova,}} \\
\small{ e-mail: \texttt{palii@jinr.ru} } }

\date{}

\maketitle

\begin{abstract}
For an adjoint action of a Lie group $G$ (or its subgroup) on Lie algebra $Lie(G)$
we suggest a method for construction of invariants. The method is easy in implementation and may shed the light on algebraical independence of invariants. The main idea is to extent automorphisms of the Cartan subalgebra to  automorphisms of the whole Lie algebra $Lie(G)$. Corresponding matrices in a linear space $V\cong Lie(G)$ define a Reynolds operator ``gathering'' invariants of torus $\mathcal{T}\subset G$ into special polynomials. A condition for a linear combination of polynomials to be $G$-invariant is equivalent to the existence of a solution for a certain system of linear equations on the coefficients in the combination.

As an example we consider the adjoint action of the Lie group $\mr{SL}(3)$ (and its subgroup $\mr{SL}(2)$) on the Lie algebra $\mf{sl}(3)$.
\end{abstract}

\newpage

\section{Introduction}\label{Sect:Introduction}

\setcounter{equation}{0}

Algorithms in invariant theory~\cite{Sturmfels,DerksenKemper} become inefficient for Lie groups interesting for physics.  To overcome this problem, we try to use the following well known fact. Every automorphism of a Lie algebra root system (in particular, an element of the Weyl group $\mathcal{W}$) defines an automorphism of the Cartan subalgebra $\mf{h}$ and, subsequently, can be extended to the automorphism of the whole Lie algebra $\mf{g}$~\cite{Humphreys}. For example, let us consider $\mf{sl}(2)$ algebra with the standard basis $X,Y,H=[X,Y]$ and define the adjoint action $\mr{ad}_A B:=[A,B]$. The operator
\[
    S=(\exp\, \mr{ad}_X)\;( \exp\, -\mr{ad}_Y)\;( \exp\, \mr{ad}_X )
\]
induces a second order inner automorphism of $\mf{sl}(2)$:
\[
    S(X)=-Y,\qquad S(Y)=-X,\qquad S(H)=-H.
\]
The action of $S$ on to the Cartan subalgebra $\mf{h}=\{H\}$ corresponds to a Weyl group $\mathcal{W}=\mathbb{Z}/2$ reflection of the 1-dimensional $\mf{sl}(2)$-root system:
\[\sigma_\alpha(\alpha)=-\alpha.\]

This paper is organized as follows. Section 2 illustrates a method on the construction of two Casimir invariants of the Lie group $SL(3)$. An adjoint representation for the Lie algebra $\mf{sl}(3)$ is done in the subsection 2.1. In 2.2, using a simple algorithm~(\cite{DerksenKemper}, \S4.3), we  construct a Hilbert basis of invariants for a maximal torus $\mathcal{T}\subset SL(3)$. Subsections 2.3 and 2.4 describe two main steps of our method. In 2.3, we construct a set of the $S$-operators extending the action of the Weyl group for the $\mf{sl}(3)$ root system and apply corresponding  Reynolds operator to the torus invariants. Derived invariants we call \textit{Weyl blocks}. In 2.4, we join them into linear combinations which are invariant relative to the adjoint action of $SL(3)$. In section 3, $SL(2)$-invariants in the algebra $\mf{sl}(3)$ are found. In 3.1, an embedding of $SL(2)$ into $SL(3)$ is used to construct all elements necessary for our method. In 3.2, we list a set of invariants and prove that it is a fundamental system of $SL(2)$-invariants in the adjoint representation of the algebra $\mf{sl}(3)$.

\section{SL(3) Invariants}

\setcounter{equation}{0}

\subsection{Linear representation of SL(3)-adjoint action on $\bs{\mf{sl}(3)}$}\label{Sect:Adjoint}

The dimension of $\mf{sl}(3)$ is equal to 8 and the rank is 2 and, correspondingly,  two $\mf{sl}(2)$-triples generate the algebra:
\begin{alignat}{3}
    & X_1=E_{12}, & \qquad   & Y_{1}=E_{21}, & \qquad & H_1=E_{11}-E_{22}, \label{generator1}\\
    & X_2=E_{23}, & \qquad   & Y_{2}=E_{32}, & \qquad & H_2=E_{22}-E_{33},\label{generator2}
\end{alignat}
where $(E_{ij})_{kl}=\delta_{ik}\delta_{jl}$ are unit matrices. These operators with two additional
\begin{equation}\label{}
    X_{3}:=[X_1,X_2], \qquad Y_{3}:=[Y_2,Y_1]
\end{equation}
form the Cartan-Weyl basis of $\mf{sl}(3)$. Below we use the following notations for the basis elements:
\begin{equation}\label{Generators}
    \begin{array}{c|c|c|c|c|c|c|c}
      Y_1 & X_1 & Y_2 & X_2 & Y_3 & X_3 & H_1 & H_2 \\
      \hline\rule{-4pt}{15pt}
      \mre_1 & \mre_2 & \mre_3 & \mre_4 & \mre_5 & \mre_6 & \mre_7 & \mre_8 \\
    \end{array}
\end{equation}
Consider $\mf{sl}(3)$ as a vector space $V$ with the basis $\{e_i\}_{i=1}^8$ corresponding to the generators $\mre_i$.
The adjoint representation of $\mf{sl}(3)$
\begin{equation}\label{ad_algebra}
    \mr{ad}_{\mre_i}:\; \mre_j\,\mapsto \, [\mre_i,\mre_j]=c_{ij}{}^k \mre_k
\end{equation}
induces a linear transformation $\mathcal{A}_i$ in $V$ by the matrix $(A_i)_{j}{}^k=-c_{ij}{}^k$:
\begin{equation}\label{LinTrA}
    \mathcal{A}_i\,(e_j):=e_k \;(\tilde{A}_i)^{k}{}_j, \qquad \tilde{A}_i=\chi^{-1} A_i \, \chi,
\end{equation}
where $k$ numerates rows, $j$ numerates columns of the matrix $(\tilde{A}_i)^{k}{}_j$ and $\chi$ is the matrix of the
Cartan metric on $\mf{sl}(3)$:
\begin{equation}\label{def: Cartan metric}
    \chi_{ij}:=\mt(A_i \cdot A_j)=c_{ip}{}^q c_{jq}{}^p.
\end{equation}
The adjoint action of the Lie group $SL(3)$ on its algebra $\mf{sl}(3)$
\begin{equation}\label{}
    \exp(t\,\mre_i) \, \mre_j \, \exp(-t\,\mre_i) = \exp(t\,\mr{ad}_{\mre_i}) \, \mre_j
    \qquad \mre_i,\mre_j \in \mf{sl}(3),
\end{equation}
can be represented as a linear transformation of $V$ by one-parameter groups:
\begin{equation}\label{AdGroupLinearTrans}
    \mathcal{C}_i\,(e_j):
    =e_k (\exp[t\,\tilde{A}_i])^{k}{}_j\,.
\end{equation}
In this section we construct polynomial $SL(3)$-invariants in variables
\begin{equation}\label{Variables}
    \begin{array}{c|c|c|c|c|c|c|c}
      e_1 & e_2 & e_3 & e_4 & e_5 & e_6 & e_7 & e_8 \\
      \hline\rule{-4pt}{15pt}
      y_1 & x_1 & y_2 & x_2 & y_3 & x_3 & h_1 & h_2 \\
    \end{array}
\end{equation}

\subsection{Torus Invariants}

A $d$-dimensional torus $\mathcal{T}\cong (\mathbb{C}^*)^d$ acting on the ring $\mathbb{C}[V]=\mathbb{C}[x_1,\ld,x_n]$ is isomorphic to a group of the diagonal matrices
\begin{equation}\label{Tor}
   T=\mr{diag}(\prod_{i=1}^d t_i^{a_{1i}},\ld,\prod_{i=1}^d t_i^{a_{ni}}),
\end{equation}
where $t_1,\ld,t_d \in \mathbb{C}^*$ are torus variables and $x_i$ have the weights
\[
    \omega_i=(a_{i,1},\ld,a_{i,d}), \qquad
    a_{ij}\in\mathbb{Z}.
\]
$\mathcal{T}$ maps monomials into monomials and its invariant ring $\mathbb{C}[V]^{\mathcal{T}}$ is monomial. Suppose $\mathcal{T}\subset \mr{G},$ for the Lie group $G$. Then the algebra of $G$-invariants is a subalgebra of torus invariants
\[ \mathbb{C}[V]^{\mr{G}}\subset \mathbb{C}[V]^{\mathcal{T}}, \]
i.e. a G-invariant polynomial is a linear combination of $\mathcal{T}$-invariant monomials.

By definition a weight $\omega_i$ of an element $X_i$ is the set of values of the corresponding root $\alpha_i$ on the elements $H_j$ of the Cartan subalgebra $\mf{h}$. The values $\alpha_i(H_j)$ can be directly read from the commutation relations in the Cartan-Weyl basis:
\begin{equation}\label{}
    [H_j,X_i]=\alpha_i(H_j)X_i.
\end{equation}

Weights of the variables~(\ref{Variables}) are ones of the generators~(\ref{Generators})

\begin{equation}\label{sl3weights}
    \begin{array}{c|c|c|c|c|c|c|c}
      y_1 & x_1 & y_2 & x_2 & y_3 & x_3 & h_1 & h_2 \\
      \hline\rule{-4pt}{15pt}
      (-2,1) & (2,-1) & (1,-2) & (-1,2) & (-1,-1) & (1,1) & (0,0) & (0,0)
    \end{array}
\end{equation}

An efficient algorithm 4.3.1 from~\cite{DerksenKemper} allows to find the \textit{Hilbert basis}~\cite{Sturmfels}, i.e. a minimal generating set, for the torus invariants. We realize it in the following way. Let us construct the $d$-dimensional set $\Omega$ of boxes marked by weights of the variables $x_i$ and 'filled in' with the corresponding $x_i$ considered as initial monomials. After that we produce new monomials degree by degree multiplying on $x_i$ monomials derived at the previous step. If a new monomial $m=x_1^{n_1}\ld x_k^{n_k}$ with the weight
\begin{equation}\label{}
    \omega_{m}=(\sum_{i=1}^k a_{i,1}n_i,\ld, \sum_{i=1}^k a_{i,d}n_i)
\end{equation}
is not divisible on derived early monomials with the same weight $\omega_m$ we put $m$ into the corresponding box. We discard a new monomial if its weight is not in $\Omega$. An output of the algorithm is the set of $\mathcal{T}$-invariant monomials with the weight $(\underbrace{0,\ld,0}_d)$.

For $\mr{T}\subset\mr{SL}(3)$ the two dimensional torus $\mr{T}$~(\ref{Tor}) is the following:
\begin{equation}\label{}
    \mr{T}=\mr{diag}\; (t_1^{-2}t_2^{1},\; t_1^{2}t_2^{-1},\; t_1^{1}t_2^{-2}\; t_1^{-1}t_2^{2},\;
    t_1^{-1}t_2^{-1},\; t_1^{1}t_2^{1},\; t_1^{0}t_2^{0},\; t_1^{0}t_2^{0}),
\end{equation}
and the Hilbert basis consists from the monomials:
\begin{equation}\label{SL3TorInvs}
    h_1,\; h_2,\; x_1y_1,\; x_2y_2,\; x_3y_3,\; x_1x_2y_3,\; x_3y_1y_2.
\end{equation}

\subsection{Weyl blocks}

The root system $\Phi_{\mf{sl}(3)}$ of $\mf{sl}(3)$ algebra lies in the Euclidean plane and has six  roots $\pm\alpha_{i},\;i=1,2,3\,$~(see figure 1). A group of the reflections $\sigma_{i}$ relative to hyperplanes orthogonal to the roots $\alpha_i$, i.e. the Weyl group $\mathcal{W}_{\mf{sl}(3)}$, is the permutation group $\mf{S}_3$~\cite{Humphreys}
\begin{equation}\label{Sigma3}
   \mf{S}_3= \{\mathbb{I},\,(12),\,(23),\,(13),\,(123),\,(132)\}
\end{equation}
The permutations $(12)$ and $(23)$ correspond to the reflections $\sigma_{1}$ and $\sigma_{2}$ (relative to hyperplanes shown by dashed lines in the figure 1) respectively and generate $\mf{S}_3$ according to formulas:
\begin{align}
    (123)&=(12)\circ(23), \label{gen123}\\
    (132)&=(23)\circ(12), \label{gen132}\\
    (13)&=(12)\circ(23)\circ(12). \label{gen13}
\end{align}

\begin{figure}
\begin{center}
\includegraphics[scale=0.6]{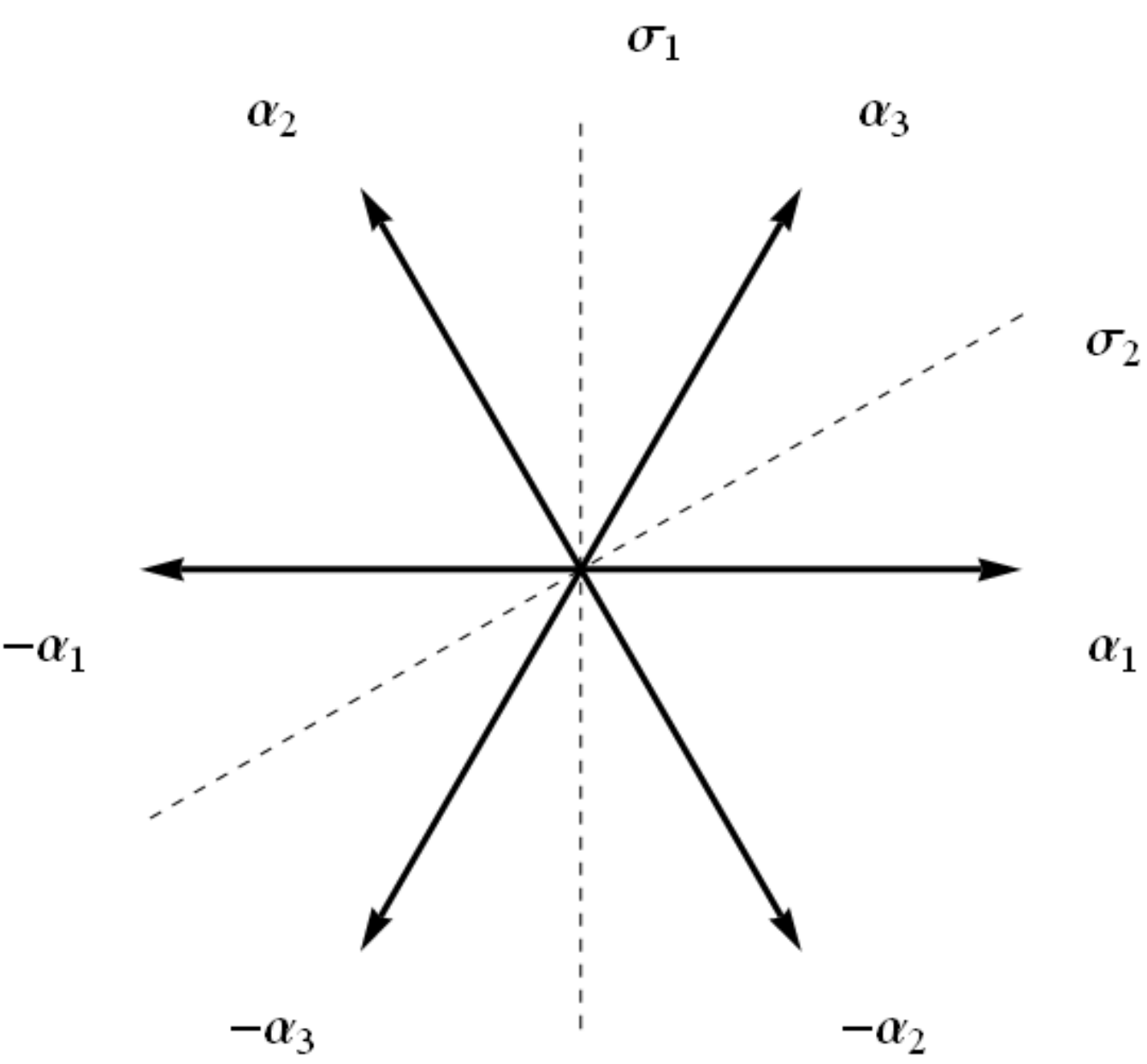}
\caption{$\mf{sl}(3)$ root system}
\end{center}
\end{figure}

Suppose we treat $\mf{sl}(3)$ algebra as a linear space $V$ according to the subsection~\ref{Sect:Adjoint}.
The reflections $\sigma_{1}$ and $\sigma_{2}$  can be extended to linear operators in $V$ analogously to $\mf{sl}(2)$ case~(see Introduction)
 \begin{align}
    & S_1=\exp[\tilde{A}_{X_1}] \exp[-\tilde{A}_{Y_1}] \exp[\tilde{A}_{X_1}] \label{s1}\\
    & S_2=\exp[\tilde{A}_{X_2}] \exp[-\tilde{A}_{Y_2}] \exp[\tilde{A}_{X_2}] \label{s2}
\end{align}
as we have two $\mf{sl}(2)$-triples~(\ref{generator1}),\,(\ref{generator2}) in $\mf{sl}(3)$. Transformations by $S_1,S_2$ of basis elements of $\mf{sl}(3)$
\begin{equation}\label{S_transformation}
    \mathcal{S}_i\,(e_j):
    =e_k (S_i)^{k}{}_j\,,
\end{equation}
are shown in the following tables:
\begin{equation}\label{S1trans} \mathcal{S}_1:\qquad\qquad
    \begin{array}{c|c|c|c|c|c|c|c}
      y_1 & x_1 & y_2 & x_2 & y_3 & x_3 & h_1 & h_2 \\ \hline
      -x_1 & -y_1 & y_3 & x_3 & -y_2 & -x_2 & -h_1 & h_3
    \end{array}
\end{equation}
\begin{equation}\label{S2trans} \mathcal{S}_2:\qquad\qquad
    \begin{array}{c|c|c|c|c|c|c|c}
      y_1 & x_1 & y_2 & x_2 & y_3 & x_3 & h_1 & h_2 \\ \hline
      -y_3 & -x_3 & -x_2 & -y_2 & y_1 & x_1 & h_3 & -h_2
    \end{array}
\end{equation}
where $h_3:=h_1+h_2$. As a proof of correctness of extension we compare the action of $\mathcal{S}_1,\mathcal{S}_2$ on the Cartan subalgebra $\mf{h}=\{h_1,h_2,h_3\}$ and Weyl group reflections of the dual roots $\alpha_i\leftrightarrow H_i,\;i=1,2,3$ in the $\mf{sl}(3)$ root space $\Phi_{\mf{sl}(3)}$:
\begin{alignat*}{2}\label{}
    & \sigma_{1}(\alpha_1)=-\alpha_1 & \qquad \mathcal{S}_{1}(h_1)&=-h_1  \\
    & \sigma_{1}(\alpha_2)=\alpha_3 & \qquad \mathcal{S}_{1}(h_2)&=h_3 \\
    & \sigma_{1}(\alpha_3)=\alpha_2 & \qquad \mathcal{S}_{1}(h_3)&=h_2
\end{alignat*}
\begin{alignat*}{2}\label{}
    & \sigma_{2}(\alpha_1)=\alpha_3 & \qquad \mathcal{S}_{2}(h_1)&=h_3  \\
    & \sigma_{2}(\alpha_2)=-\alpha_2 & \qquad \mathcal{S}_{2}(h_2)&=-h_2  \\
    & \sigma_{2}(\alpha_3)=\alpha_1 & \qquad \mathcal{S}_{2}(h_3)&=h_1
\end{alignat*}
Taking a Cartan decomposition
\begin{equation}\label{}
    \mf{sl}(3)=\{h_1,h_2\}\oplus\coprod_{\alpha\in\Phi_{\mf{sl}(3)}}L_\alpha
\end{equation}
we see from the tables~(\ref{S1trans}),(\ref{S2trans}) that every one-dimensional subspace $L_\alpha$ is transformed according to the reflection $\sigma_i$: $\;\mathcal{S}_i:\;L_{\alpha}\rightarrow L_{\sigma_i(\alpha)}$.

Let us consider matrices
\begin{alignat}{2}\label{}
    & S_0=\mathbb{I}, & \qquad & S_3=S_1 S_2 S_1, \\
    & S_1, & \qquad & S_4=S_1 S_2, \\
    & S_2, & \qquad & S_5=S_2 S_1,
\end{alignat}
where we construct $S_3,S_4,S_5$ following to~(\ref{gen13}),(\ref{gen123}),(\ref{gen132}) respectively.
$S_0,\ld,S_5$ are not a representation of $\mathcal{W}_{\mf{sl}(3)}=\mf{S}_3$~\cite{Humphreys}. $S_1,S_2$ generate a non-abelian matrix group of order 24 but we postpone investigation of this group.

We define a Reynolds operator $Rey$ on a polynomial $P$ in variables~(\ref{Variables})
\begin{equation}\label{Rey}
    \mathit{Rey}(P)=\sum_{S} \mathcal{S}_i(P), \qquad
    S=\{\mathcal{S}_0,\,\mathcal{S}_1,\,\ld,\,\mathcal{S}_5\}
\end{equation}
where transformations $\mathcal{S}_i$ correspond to the matrices $S_i$ according to~(\ref{S_transformation}).

Application of $\mathit{Rey}$ to the torus invariants~(\ref{SL3TorInvs}) yields polynomials which we call \textit{Weyl blocks} $w_{i,j}$ ($j$ numerates Weyl blocks of degree $i$). Linear combinations of the torus invariants give \textit{initial} Weyl blocks:
\begin{equation}\label{InitWB}
    w_{2,1}=x_1y_1+x_2y_2+x_3y_3, \qquad w_{3,1}=x_1x_2y_3+x_3y_1y_2.
\end{equation}
Taking products of torus invariants we derive additional blocks:
\begin{align}\label{ProdWB}
    & w_{2,2}=h_1^2+h_1 h_2 +h_2^2, \\
    & w_{3,2}=2h_1^3+3h_1^2 h_2 - 3h_1 h_2^2 - 2h_2^3, \\
    & w_{3,3}=x_1y_1(h_1+2h_2)-x_2y_2(2h_1+h_2)+x_3y_3(h_1-h_2).
\end{align}

\subsection{$\bs{\mr{SL}(3)}$ Invariants}

The group $\mr{SL}(3)$ is generated by eight one-parameter subgroups $\mathcal{C}^{(k)},\;k=1,\ld,8$~(\ref{AdGroupLinearTrans}) corresponding to the basis elements~(\ref{Generators}). A difference of a Weyl block $w_{i,j}$ under the action of $\mathcal{C}^{(k)}$  has a part which is linear in the parameter $t$:
\begin{equation}\label{}
    \Delta^{(k)}_{i,j}:=\frac{d}{dt}\mathcal{C}^{(k)}(w_{i,j})\Big|_{t=0}
\end{equation}
A homogeneous linear combination $\sum_{i,j} \mu_{i,j} w_{i,j}\;$ of the Weyl blocks with numerical coefficients $\mu_{i,j}$ is $\mr{SL}(3)$-invariant if and only if the sum $\sum_{i,j} \mu_{i,j} \Delta^{(k)}_{i,j}\;$ is equal to zero for all $\mathcal{C}^{(k)}$.

Now let us take into account that $\Delta^{(k)}_{i,j}$ are polynomials in the variables $e_i,\;i=1,\ld,8$~(\ref{Variables}),
\begin{equation*}
    \Delta^{(k)}_{i,j}(e_1,\ld,e_8)=\sum_{\mr{d}\in \mr{Deg}} \phi^{(k)}_{i,j}(\mr{d})\,e_1^{d_1}\ld e_8^{d_8}, \qquad
    \mr{d}=(d_1,\ld,d_8).
\end{equation*}
where $\mr{Deg}$ is a set of multidegrees $\mr{d}$ of monomials $e_1^{d_1}\ld e_8^{d_8}$ in $\Delta^{(k)}_{i,j}$. Rewriting the equation
\begin{equation}\label{}
    \sum_{i,j} \mu_{i,j} \Delta^{(k)}_{i,j}=0
\end{equation}
as a sum over $\mr{d}$, we get a homogeneous system of linear equations (one equation for every monomial $e^{\mr{d}}$ and fixed $k$)
\begin{equation}\label{SystemMu}
    \sum_{i,j} \mu_{i,j} \phi^{(k)}_{i,j}(\mr{d})=0
\end{equation}
for the coefficients $\mu_{i,j}$.

For the Weyl blocks~(\ref{InitWB}),~(\ref{ProdWB}) we get equations:

second order $w_{2,i}$,
    \begin{equation}\label{}
    \mu_{2,1}-3\mu_{2,2}=0,
    \end{equation}

third order $w_{3,i}$,
    \begin{equation}\label{}
    \mu_{3,1}-3\mu_{3,3}=0,\qquad
    9\mu_{3,2}-\mu_{3,3}=0.
    \end{equation}

The equations have solutions spanned by vectors
\begin{equation}\label{nullvectros}
    \{\mu_{2,1},\mu_{2,2}\}=\{3,1\}, \qquad \qquad
    \{\mu_{3,1},\mu_{3,2},\mu_{3,3}\}=\{27,1,9\}
\end{equation}

Thus we derive two Casimir invariants in the algebra $\mf{sl}(3)$ under the adjoint action of the group $\mr{SL}(3)$:
\begin{align}
    & \mf{C}_2=3 w_{2,1} + w_{2,2}, \label{SL3Cas2}\\
    & \mf{C}_3=27 w_{3,1} + w_{3,2} + 9 w_{3,3}. \label{SL3Cas3}
\end{align}

The Weyl blocks $w_{2,2}$ and $w_{3,2}$ represent reductions of the Casimir invariants onto the Cartan subalgebra $\mf{H}$ because the variables $h_1,h_2$ correspond to the generators $H_1,H_2\in\mf{H}$. We substitute  $h_i$ by dual roots $\alpha_i$ and get invariants of the Weyl group
\begin{align}\label{}
    & w_{2,2}\;\rightarrow \mathcal{I}^{\mathcal{W}}_2=
    \alpha_1^2+\alpha_1 \alpha_2 +\alpha_2^2, \\
    & w_{3,2}\;\rightarrow \mathcal{I}^{\mathcal{W}}_3=
    2\alpha_1^3+3\alpha_1^2 \alpha_2 - 3\alpha_1 \alpha_2^2 - 2\alpha_2^3
\end{align}
in accordance with the Chevalley theorem about the homomorphism of Lie group invariants and Weyl group invariants~\cite{Humphreys}.

\section{SL(2) Invariants in $\bs{\mf{sl}(3)}$ algebra}

\setcounter{equation}{0}

\subsection{$\bs{\mf{sl}(2)}$ as a subalgebra of $\bs{\mf{sl}(3)}$}\label{sl2_in_sl3}

The generators $\{Y_1,H_1,X_1\}$~(see~(\ref{generator1})) form the subalgebra $\mf{su}(2)$ in $\mf{su}(3)$. We decompose $\mf{sl}(3)$ into a direct sum of irreducible $\mf{sl}(2)$-modules looking at the commutation relations:
\begin{alignat}{3}
    & [H_1,Y_2]=Y_2 & \qquad & [Y_1,Y_2]=-Y_3 & \qquad & [X_1,Y_2]=0, \\
    & [H_1,Y_3]=-Y_3 & \qquad & [Y_1,Y_3]=0 & \qquad & [X_1,Y_3]=-Y_2.
\end{alignat}
\begin{alignat}{3}\label{}
    & [H_1,X_3]=X_3 & \qquad & [Y_1,X_3]=X_2 & \qquad & [X_1,X_3]=0, \\
    & [H_1,X_2]=-X_2 & \qquad & [Y_1,X_2]=0 & \qquad & [X_1,X_2]=X_3.
\end{alignat}
It is convenient to introduce a new operator $H_0$,
\begin{equation}\label{}
    H_0:=\tfrac{1}{2}H_1+H_2,
\end{equation}
with commutators:
\begin{equation}\label{}
    [H_1,H_0]=[Y_1,H_0]=[X_1,H_0]=0.
\end{equation}
According to the above relations, $\mf{su}(3)$ is a direct sum of four $\mf{sl}(2)$-modules:
\begin{equation}\label{}
    V_0 \oplus V'_1 \oplus V''_1 \oplus V_2
\end{equation}
where each module $V_\lambda$ is a disjoint union of one-dimensional subspaces
\begin{equation}\label{}
    V_\lambda=\coprod_{\mu}V_{\lambda,\mu}, \qquad \mu=\lambda,\lambda-2,\ld,-(\lambda-2),-\lambda\;.
\end{equation}
The decomposition can be summarized by the following table.

\bigskip

\begin{equation}\label{sl3_decomposition}
\begin{array}{c|c|c|c|c}
  \mu \diagdown V_\lambda  & \quad V_2 \quad & \quad V'_1 \quad & \quad V''_1 \quad & \quad V_0 \quad\\[2mm]  \hline\rule{-4pt}{15pt}
  2 & X_1 &  &  &  \\[2mm] \hline \rule{-4pt}{15pt}
  1 &  & Y_2 & X_3 &  \\[2mm] \hline \rule{-4pt}{15pt}
  0 & H_1 &  &  & H_0 \\[2mm] \hline  \rule{-4pt}{15pt}
  -1 &  & Y_3 & X_2 &  \\[2mm] \hline  \rule{-4pt}{15pt}
  -2 & Y_1 &  &  &
\end{array}
\end{equation}

\bigskip

Embedding $\mf{sl}(2)$ into $\mf{sl}(3)$ allows us to use many constructions from the previous section. Let us reorder $\mf{su}(3)$-generators~(\ref{Generators}) in accordance with the table~(\ref{sl3_decomposition})
\begin{equation}\label{SL3moduleGen}
    Y_1,\,X_1,\,H_1,\,Y_2,\,Y_3,\,X_2,\,X_3,\,H_0
\end{equation}
The group $\mr{SL}(2)$ is generated by three one-parameter subgroups \[\mathcal{C}_{Y_1},\,\mathcal{C}_{X_1},\,\mathcal{C}_{H_1}\]
(see (\ref{AdGroupLinearTrans})) and acts on the row
\begin{equation}\label{vector}
    v=\{y_1,x_1,h_1,y_2,y_3,x_2,x_3,h_0\}
\end{equation}
from the right.

The Cartan subalgebra of $\mf{sl}(2)$ is one-dimensional: $\mf{H}_{\mf{sl}(2)}=\{H_1\}$. And the weights of the basis elements~(\ref{SL3moduleGen}) are given by the first component of their $\mf{sl}(3)$-weights~(\ref{sl3weights}) except $H_0$ with the weight $0$. From other hand, $\mf{sl}(2)$-weights are given by the value of $\mu$ in the table(\ref{sl3_decomposition}). Accordingly, one-dimensional torus $\mathcal{T}_{\mr{SL}(2)}$
\begin{equation}\label{}
    \mathcal{T}_{\mr{SL}(2)}=\mr{diag}(t^{-2},\,t^2,\,t^0,\,t^1,\,t^{-1},\,t^{-1},\,t^1,\,t^0)
\end{equation}
acts on 8-dimensional vector $v$~(\ref{vector}) from the right. The Hilbert basis of the torus invariants includes new monomials
\begin{equation}\label{SL2TorInvs}
    \begin{array}{llll}
      y_2y_3, & x_2x_3, & & \\
      y_1y_2^2, & y_1x_3^2, & x_1x_2^2, & x_1y_3^2
    \end{array}
\end{equation}
in addition to the invariants of $\mathcal{T}_{\mr{SL}(3)}$~(\ref{SL3TorInvs}) where $h_0$ must be included instead of $h_2$.

The root system $\Phi_{\mf{sl}(2)}$ is represented by $\pm\alpha_1$ and has the Weyl group $\mathcal{W}_{\mf{sl}(2)}=\mathbb{Z}_2$. The operator $s_1$~(\ref{s1}) has the block-diagonal form
\begin{equation}\label{}
    \left(
      \begin{array}{cccccccc}
        0 & -1 & 0 & 0 & 0 & 0 & 0 & 0 \\
        -1 &  & 0 & 0 & 0 & 0 & 0 & 0 \\
        0 & 0 & -1 & 0 & 0 & 0 & 0 & 0 \\
        0 & 0 & 0 & 0 & -1 & 0 & 0 & 0 \\
        0 & 0 & 0 & 1 & 0 & 0 & 0 & 0 \\
        0 & 0 & 0 & 0 & 0 & 0 & -1 & 0 \\
        0 & 0 & 0 & 0 & 0 & 1 & 0 & 0 \\
        0 & 0 & 0 & 0 & 0 & 0 & 0 & 1 \\
      \end{array}
    \right)
\end{equation}
Now the Reynolds operator~(\ref{Rey}) includes only summation over $s_1$ and the identity operator.

\subsection{$\bs{\mathrm{SL}(2)}$-Invariant polynomials in $\bs{\mf{sl}(3)}$ algebra}

Following the described above procedure, we get as a result:
\begin{equation}\label{SL2invariants}
    \begin{array}{c|l}
    \mbox{degree} & \mbox{invariants} \\[2mm] \hline\rule{0pt}{15pt}
  1 & \mathcal{I}_1=\underline{h_0} \\[2mm] \hline\rule{0pt}{15pt}
  2 & \mathcal{I}_2=h_1^2+4\,\underline{x_1 y_1} \\[2mm]
    & \mathcal{I}_3=\underline{x_2 y_2 + x_3 y_3} \\[2mm] \hline\rule{0pt}{15pt}
  3 & \mathcal{I}_4=h_1 y_2 y_3+\underline{y_1 y_2^2-x_1 y_3^2} \\[2mm]
    & \mathcal{I}_5=h_1 x_2 x_3+\underline{x_1 x_2^2-x_3^2 y_1} \\[2mm]
    & \mathcal{I}_6=h_1(x_2 y_2 - x_3 y_3)-2\,(\underline{y_1y_2x_3+x_1x_2y_3})
\end{array}
\end{equation}
where initial Weyl blocks are underlined.

We shall prove that the set of invariants $\mathcal{I}_i,\;i=1,\ld,6$ generates the invariant ring $\mathbb{C}[y_1,x_1,\ld,h_0]^{\mr{SL}(2)}$, i.e. this set is a fundamental system of $\mr{SL}(2)$-invariants in $\mf{sl}(3)$ algebra. This is the case if and only if the Hilbert series $H(R,q)$ of the subalgebra $R=\mathbb{C}[\mathcal{I}_1,\ld,\mathcal{I}_6]$ is equal to the Molien series~\cite{Sturmfels} of $\mathbb{C}[y_1,x_1,\ld,h_0]^{\mr{SL}(2)}$.

For this purpose we calculate the Molien function~\cite{DerksenKemper} for the adjoint action of $\mr{SU}(2)$ (i.e. the real compact form of $\mr{SL}(2)$) on the algebra $\mf{su}(3)$. An injective map of $\mathrm{SU}(2)$ into $\mathrm{SU}(3)$ \begin{equation}\label{}
  \rho=\begin{pmatrix}
          \mathrm{SU}(2) & 0 \\
          0 & 1 \\
        \end{pmatrix}\;\subset \;\mathrm{SU}(3)
\end{equation}
corresponds to our choice of representation of $\mf{sl}(2)$ in $\mf{sl}(3)$ (see the previous subsection~\ref{sl2_in_sl3}). One can instead of the adjoint action
\begin{equation}\label{}
    \rho\,a\,\rho^{-1} \qquad a\in\mf{su}(3)
\end{equation}
consider a linear representation $L$ of $\rho$ on $V=\mathbb{R}^{9}$
\[ V_A^\prime = L_{AB}V_{B}\,,
\qquad\qquad  L_{AB} \in
\pi(\rho)\otimes\overline{\pi(\rho)} \,, \]
where $\overline{\pi(\rho)}$ is the complex conjugation of $\pi(\rho)$. After that we calculate the Molien function
\begin{equation}\label{eq:Molienfunct}
    M(\mathbb{C}[V]^{\mathrm{SU}(2)},q)=\int_{\mathrm{SU}(2)}\,
    \frac{d\mu}{\det(\mathbb{I}-q\,\pi(\rho)\otimes\overline{\pi(\rho)})}, \qquad
|q|<1\,,
\end{equation}
where  ${d\mu}$ is a Haar measure for $\mathrm{SU}(2)$, $\mathbb{I}$ is an identity operator, and $\pi(\rho)$ is a representation of $\rho$. Using Weyl integration formula~\cite{Weyl}, one can reduce the integral over $\mathrm{SU}(2)$ to the integral over its maximal torus $\mathcal{T}_{\mathrm{SU}(2)}$. An expression of ${d\mu}$ via the Haar measure $d\phi$ on $\mathcal{T}_{\mathrm{SU}(2)}$ is the following one:
\begin{equation}\label{Haar}
    \begin{split}
      d\mu= & \frac{1}{\pi}\sin^2(\phi)\,d\phi  \qquad 0\leq \phi \leq 2\pi\\
       = & -\frac{1}{4\pi i}\frac{dz}{z^3}(1-z^2)^2 \qquad \mbox{where} \quad z=e^{i\phi}.
    \end{split}
\end{equation}
$\pi(\rho)$ has a diagonal form
\begin{equation}\label{}
    \pi(\rho)=
        \begin{pmatrix}
          z & 0 &  0 \\
          0 & z^{-1}& 0 \\
          0 & 0 & 1 \\
        \end{pmatrix}\;,
\end{equation}
and, hence,
\begin{equation}\label{}
    \begin{split}
      \pi(\rho)\otimes\overline{\pi(\rho)}= & \mr{diag}(z,z^{-1},1)\otimes(z^{-1},z,1) \\
        = & \mr{diag}(1,z^2,z,z^{-2},1,z^{-1},z^{-1},z,1),
    \end{split}
\end{equation}
\begin{equation}\label{Det}
    \det||\mathbb{I}-q\pi_{{}_G}\otimes \bar{\pi}_{{}_G}||=
    (1-q)^3(1-qz)^2(1-qz^2)(1-qz^{-1})^2(1-qz^{-2}).
\end{equation}
Substituting (\ref{Haar}), (\ref{Det}) in the expression (\ref{eq:Molienfunct}) we derive that the Molien function is equal to the following contour integral over the unit circle:
\begin{equation}\label{}
    M(\mathbb{C}[V]^{\mathrm{SU}(2)},q)= -\frac{1}{4\pi i}\frac{1}{(1-q)^3}
\oint \frac{zdz(1-z^2)^2}{(1-qz)^2(1-qz^2)(z-q)^2(z^2-q)}.
\end{equation}
There are two simple poles $z=\pm\sqrt{q}$ and a double pole $z=q$ (of order $m=2$). Application of the residue theorem~\cite{DerksenKemper} yields the result:
\begin{equation}\label{Molien_result}
    M(\mathbb{C}[V]^{\mathrm{SU}(2)},q)=\frac{1+q^3}{(1-q)^2(1-q^2)^2(1-q^3)^2}.
\end{equation}
$M(\mathbb{C}[V]^{\mathrm{SU}(2)},q)$ satisfies a functional equation
\begin{equation}\label{}
    M(\mathbb{C}[V]^{\mathrm{SU}(2)},q^{-1})=q^{9}  M(\mathbb{C}[V]^{\mathrm{SU}(2)},q) \,,
\end{equation}
hence,  $\mathbb{C}[V]^{\mathrm{SU}(2)}$ being a graded Cohen-Macaulay algebra~\cite{HochsterRoberts} is the Gorenstein one according to Stanley~\cite{Stanley}.

Required Molien series of $\mathrm{SU}(2)$-invariant ring $\mathbb{C}[V]^{\mathrm{SU}(2)}$ is nothing then
an expansion of $M(\mathbb{C}[V]^{\mathrm{SU}(2)},q)$ in series of powers of $q$. Now we turn out to the calculation of the Hilbert series $H(R,q)$ for the subalgebra $R=\mathbb{C}[\mathcal{I}_1,\ld,\mathcal{I}_6]$. Using slack-variables method~\cite{Sturmfels}, we find one syzygy between the invariants $\mathcal{I}_i$~(\ref{SL2invariants}),
\begin{equation}\label{syzygy}
    \mathcal{I}_2 \mathcal{I}_3^2-4 \mathcal{I}_4 \mathcal{I}_5-\mathcal{I}_6^2\,
\end{equation}
thus $\{\mathcal{I}_1,\ld,\mathcal{I}_5\}$ are algebraically independent while $\mathcal{I}_6$ can be considered as algebraically dependent. This implies that $R$ is decomposed as the direct sum of graded $\mathbb{C}$-vector space
\begin{equation}\label{Hironaka_decomposition}
    \mathbb{C}[\mathcal{I}_1,\ld,\mathcal{I}_6]=
    \mathbb{C}[\mathcal{I}_1,\ld,\mathcal{I}_5]\oplus
    \mathcal{I}_6\,\mathbb{C}[\mathcal{I}_1,\ld,\mathcal{I}_5]
\end{equation}
with the Hilbert series
\begin{equation}\label{Hilbert series Invariants}
    H(R,q)=\frac{1+q^3}{(1-q)(1-q^2)^2(1-q^3)^2}.
\end{equation}
We see that the denominator of the Molien function~(\ref{Molien_result}) has an extra term $1-q$ in comparison with $H(R,q)$. The reason is that we have calculated $M(\mathbb{C}[V]^{\mathrm{SU}(2)},q)$ in 9-dimensional vector space but $\mf{sl}(3)$ has dimension 8. This completes the proof.

As a corollary we get that a Hironaka decomposition~\cite{Sturmfels} of $R$ is done by the formula~(\ref{Hironaka_decomposition}).  $\mathcal{I}_1,\ld,\mathcal{I}_5$ are primary invariants and $\mathcal{I}_6$ is a secondary invariant. For example, a decomposition of the Casimir invariants~(\ref{SL3Cas2}) is the following:
\begin{align}\label{}
    \mf{C}_2&=\mathcal{I}_1^2+\tfrac{3}{4}\, \mathcal{I}_2+3 \,\mathcal{I}_3,  \\
    \mf{C}_3&=-2\, \mathcal{I}_1^3+\tfrac{9}{2}\,\mathcal{I}_1 \mathcal{I}_2
    -9\, \mathcal{I}_1 \mathcal{I}_3-\tfrac{27}{2}\,\mathcal{I}_6.
\end{align}

\section{Summary}

For construction of Lie group invariants we propose a method based on the Weyl block structure of invariants. To produce $\mathcal{W}$-blocks, we extend the Weyl group action on the Cartan subalgebra to the automorphism of the whole Lie algebra. Corresponding set of operators defines a Reynolds operator $Rey$. Being applied to Hilbert basis of torus invariants operator $Rey$ produces initial $\mathcal{W}$-blocks. Beside that we construct $\mathcal{W}$-blocks from all possible products of torus invariants up to the higher degree of invariants from the Hilbert basis. A homogeneous linear combination of $\mathcal{W}$-blocks is invariant relative to the adjoint action of a Lie group if a difference of the combination under a transformation is zero. This condition gives a system of linear equations. If there is a solution of the system, then we get an invariant.

Using the described above method, we construct fundamental sets of invariants of the adjoint actions of the Lie groups $\mr{SL}(3)$ and $\mr{SL}(2)$ on the Lie algebra $\mf{sl}(3)$. We see that the Weyl block structure of invariants reflects very deep properties of Lie algebras and Lie groups such as generating relations, properties of root systems, irreducible representations and so on. All this information can be useful for studying rings of invariants. For example, let us mention that the initial $\mathcal{W}$-blocks of the invariants~(\ref{SL2invariants}) taken with their coefficients satisfy the same syzygy~(\ref{syzygy}) as invariants themselves.

\section*{Acknowledgments}

The author thanks Vladimir Gerdt and Arsen Khvedelidze for helpful discussions.
The work was supported in part by the RFBR  (grant No. 10-01-00200)
and by the Ministry of Education and Science of the Russian Federation
(grant No. 3802.2012.2).

%
%
%
%



\begin{thebibliography}{99}

\bibitem{Sturmfels} B.~Sturmfels, \textit{Algorithms in Invariant Theory}.
2nd edition, Springer-Verlag, 2008.

\bibitem{DerksenKemper} H.~Derksen and G.~Kemper,
\textit{Computational Invariant Theory}, Encyclopedia of
Mathematical Sciences, vol. 130, Springer-Verlag, Berlin, 2002.

\bibitem{Humphreys} J.E.~Humphreys, \textit{Introduction to Lie Algebras and
Representation Theory}. Springer-Verlag, N.-Y., Inc., 1978.

\bibitem{HochsterRoberts} M. Hochster and J. Roberts, Rings of invariants of reductive groups acting on regular rings are Cohen-Macaulay, Advances in Mathematics 13 (1974) pp 125-175.

\bibitem{Weyl} H.~Weyl, \textit{The classical groups - their invariants and representations}. Princeton University press, 1939. Princeton

\bibitem{Stanley} R. Stanley, Invariants of finite groups and their applications to combinatorics, Bulletin of the AMS, v1 No 3 (1979) pp. 475-511.

\end{thebibliography}
\end{document}